\documentclass[12pt,letter]{amsart}
\usepackage{amsfonts,amssymb,amsthm,amsmath,amsxtra,amscd}
\usepackage[all]{xy}
\usepackage[dvips]{graphics}

\theoremstyle{plain}
\newtheorem{thm}{Theorem}[subsection]
\newtheorem{lem}[thm]{Lemma}
\newtheorem{prop}[thm]{Proposition}
\newtheorem{cor}[thm]{Corollary} 
  
\newtheorem*{thmA}{Theorem A}
\newtheorem*{thmB}{Theorem B}
\newtheorem*{propC}{Proposition C}
\newtheorem*{propD}{Proposition D}
\newtheorem*{thmE}{Theorem E}
\newtheorem*{thmF}{Theorem F}

\theoremstyle{definition}
\newtheorem{defi}[thm]{Definition}
\newtheorem{nt}[thm]{Notation}

\theoremstyle{remark}
\newtheorem{rmk}[thm]{Remark}
\newtheorem{eg}[thm]{Example}
\newtheorem*{egE1}{Example {\rm E1}}
\newtheorem*{egE2}{Example {\rm E2}}
\newtheorem*{egE3}{Example {\rm E3}}
\newtheorem*{egE4}{Example {\rm E4}}

\newcommand{\R}{{\mathbb{R}}}                            
\newcommand{\C}{{\mathbb{C}}}

\newcommand{\Z}{{\mathbb{Z}}}

\renewcommand{\P}{{\mathbb{P}}}

\newcommand{\F}{{\mathbb{F}}}

\renewcommand{\O}{\mathcal O}

\newcommand{\LL}{\mathcal L}

\renewcommand{\a} {\alpha}
\renewcommand{\b} {\beta}
\newcommand{\g} {\gamma}

\newcommand{\f} {\phi}
\newcommand{\ff} {\psi}
\newcommand{\e} {\eta}
                       
\renewcommand{\l} {\lambda}
\newcommand{\n} {\nu}

\newcommand{\p} {\pi}
\renewcommand{\r} {\rho}
\newcommand{\s} {\sigma}
\renewcommand{\t} {\tau}

\newcommand{\G} {\Gamma}

\renewcommand{\S} {\Sigma}

\renewcommand{\.} {\cdot}
\newcommand{\ov} {\overline}

\renewcommand{\~} {\widetilde}

\newcommand{\rat}{\dashrightarrow}                      %arrows

\DeclareMathOperator{\PGL}  {PGL}                     
\DeclareMathOperator{\Gr} {Gr}
\DeclareMathOperator{\NE} {NE}
\DeclareMathOperator{\CNE} {\ov{\NE}}

\DeclareMathOperator{\Spec} {Spec}
\DeclareMathOperator{\Pic} {Pic}
\DeclareMathOperator{\Aut} {Aut}
\DeclareMathOperator{\Bir} {Bir}
\DeclareMathOperator{\Sing} {Sing}
\DeclareMathOperator{\Bl} {Bl}
\DeclareMathOperator{\NFC} {NFC}

\DeclareMathOperator{\T} {T}
\DeclareMathOperator{\NS} {NS}
\DeclareMathOperator{\cont} {cont}
\DeclareMathOperator{\Proj} {Proj}
 
\title [On planar Cremona maps]
{On planar Cremona maps of prime order}
\author[T. de Fernex]{Tommaso de Fernex}
\address{Department of Mathematics, University of Michigan,
East Hall, 525 East University Avenue, Ann Arbor, MI 48109-1109, USA}
\email{defernex@umich.edu}
\address{Fax: +1 (734) 763-0937}

\begin{document}

\maketitle

\begin{abstract}
This paper contains a new proof of the classification
of prime order elements of $\Bir(\P^2)$ up to conjugation.
The first results on this topic can be traced back to classic works by 
Bertini and Kantor, among others. The innovation introduced by this paper
consists of explicit geometric constructions of these Cremona 
transformations and the parameterization
of their conjugacy classes. The methods employed here are inspired
to~\cite{BB:99}, and rely on the
reduction of the problem to classifying prime order automorphisms
of rational surfaces.
This classification is completed by combining equivariant
Mori theory to the analysis of the action on
anticanonical rings, which leads to characterize
the cases that occur by explicit equations (see~\cite{Zhang:00}
for a different approach). 
Analogous constructions in higher dimensions are also discussed.
\end{abstract}

\section*{Introduction}

One of the first contributions to the classification of conjugacy
classes in the Cremona group of $\P^2$ can be attributed to Bertini
for his work on birational involutions~\cite{Bertini:1877}.
The classification of all finite subgroups of $\Bir(\P^2)$ up to conjugation
was successively completed by Kantor in~\cite{K:1895}.
On the same topic, one 
also finds~\cite{Wiman:1897},~\cite{Autonne:1885},~\cite{Godeaux:49} 
and~\cite{Engel:56}.
The classification of finite order planar Cremona maps,
up to conjugation, is equivalent to the classification of 
normal multiple rational planes, up to birational equivalence. 
In this area one can find the results of Bottari~\cite{Bottari:99} and
Castel\-nuovo and Enriques~\cite{CE:00}.
Recently, Bayle and Beauville~\cite{BB:99} and 
Calabri~\cite{Calabri:00},~\cite{Calabri:01} gave 
new proofs of the birational classification of, respectively,
involutions in $\Bir(\P^2)$ and double and triple rational planes.
Closely related to these topics are 
results leading towards the determination of
automorphism groups of rational surfaces, such 
as~\cite{Segre:42},~\cite{Manin:74}, \cite{Manin:67}, \cite{Gizatullin:80},
\cite{Koitabashi:88}, \cite{Hosoh:96},~\cite{Hosoh:97}~\cite{Zhang:00} 
and~\cite{Zhang:02}.
We refer to~\cite{Alberich:01} for an account of
the classic theory of planar Cremona maps.

The purpose of this paper is to give a new proof of the birational
classification of planar Cremona maps of prime order,
aiming a better understanding of the geometry 
governing these transformations. 
Elements representing each conjugacy class will be
constructed by first realizing them as automorphisms
on birational equivalent models, and then interpreting 
the constructions in terms of the geometry of $\P^2$.
Parameterization of their classes will follow from 
this approach.

Extending the methods in~\cite{BB:99}, 
the classification of Cremona maps of prime order
is reduced, through a suitable resolution of indeterminacy,  
to that of automorphisms of prime order of smooth rational surfaces. 
In fact, we will deal with automorphisms
of prime order of smooth projective surfaces whose canonical
class is not numerically effective.
The classification splits into two categories, according to the 
rank of the invariant part of the Ner\'on-Severi group of the surface.
If this rank is at least 2, we apply Mori theory in the spirit 
of~\cite{BB:99} and~\cite{Zhang:00}, 
searching for equivariant fibration structures.
Otherwise the rank is 1. Then, after observing that the surface is Del Pezzo,
the classification is completed by considering
the action that the automorphism induces on the anticanonical ring.
This approach enables us to characterize each case 
by explicit equations and identify families of analogous
automorphisms in all dimensions. 

The classification of automorphisms of prime order of smooth rational surfaces
has been already proved, by different methods, 
by Dolgachev and Zhang in their very nice paper~\cite{Zhang:00}. 
We would like to mention that 
Theorem~A below differs from~\cite[Theorem~1]{Zhang:00}
in the way certain cases are characterized: in~\cite{Zhang:00} 
surfaces and automorphisms are constructed
and characterized as cyclic coverings over their quotients,
whereas in this paper we characterize
them in terms of their equations.

This paper is organized as follows. The main results
of classification, given 
in Theorems~A,~B,~E,~F, are stated in Section~1.
Sections~2 and~3 are respectively devoted to fix the notation and
present some preliminary material. Section~4 contains the proof of 
Theorems~A and ~B. In Sections~5 and~6, we see two more properties
concerning automorphisms of surfaces: Propositions~C and~D.
Finally, in Section~6, Theorems~A,~B and 
Propositions~C,~D are applied to prove Theorems~E and~F.
Special numeration, labeling certain cases, will be consistently
adopted in all statements.

\subsection*{Acknowledgments}

First and foremost, I would like to express my gratitude 
to Professor L.~Ein for his invaluable advice and help throughout 
the preparation of this paper. 
It is a pleasure to thank A.~Beauville, N.~Budur, A.~Lanteri
and R.~Lazarsfeld for several useful discussions, and
I.~Dolgachev and D.-Q.~Zhang
for precious comments and for sending me their preprint while I was
working on this problem. 
I would like to thank the referee for his corrections
and for indicating an alternative argument that simplifies the 
original proofs of Proposition~\ref{d=3} and Lemma~\ref{trivial}.
I am grateful to A.~Albano, M.~Alberich, A.~Calabri, 
L.~Caporaso, L.~Chiantini, F.~Russo, and F.L.~Zak 
for very interesting comments.
I would like to thank the University of Genova, the
University of Hong Kong and the University of Illinois
at Chicago for their hospitality, and 
the MURST of the Italian Government
for partial support in the framework of
the National Research Project (Cofin 2000)
``Geometry on Algebraic Varieties''.

\section{The main results}

\subsection{Automorphisms of prime order of surfaces}

We work over the field of complex numbers. 
Let $X$ be a smooth projective surface, and 
$\s \in \Aut(X)$. The pair $(X,\s)$ is said to be {\it minimal} if for any 
birational morphism $\f : X \to X'$ such that $X'$ is smooth
and $\f \s \f^{-1} \in \Aut(X')$, $\f$ is an isomorphism.
Examples of minimal pairs are given by the following two 
classic involutions, whose constructions we recall here 
for the convenience of the reader. If $X$ is a smooth
Del Pezzo surface of degree 2, the linear system $|-K_X|$ defines
a double covering over $\P^2$,
branched along a smooth quartic curve; the involution 
defined by this cover is called {\it Geiser involution}.
Similarly, if $X$ is a smooth Del Pezzo surface of degree 1,
the linear system $|-2K_X|$ defines
a double covering over a quadric cone,
branched along the vertex of the
cone and a smooth curve of genus 4, and the 
corresponding involution is the {\it Bertini involution} of $X$.

\begin{thmA}
Let $X$ be a smooth projective surface whose canonical 
class is not nef, and $\s \in \Aut(X)$ be an element 
of prime order $n$ such that the pair $(X, \s)$
is minimal. Then either $(X,\s)$ is one of the following
(where any value of $n$ may occur):
\begin{enumerate}
\item[1.]   $X \cong \P^2$ and $\sigma \in \PGL(3)$;
\item[2.]   $X$ is a geometrically ruled surface and $\s$ is fiberwise,
	either inducing an effective automorphism
        on the base curve of the ruling or restricting to an effective
        automorphism on each fiber;
\end{enumerate}
or $n=2$ and $(X,\s)$ is one of the following:
\begin{enumerate}
\item[3.]   $X$ is a conic bundle and $\s$ restricts to an effective
        involution on each fiber; the two components of each singular
        fiber are flipped by $\s$;
\item[4.]   $X \cong \P^1 \times \P^1$ and $\s$ is the involution
        swapping the two rulings of $X$;
\item[5.]   $X$ is a Del Pezzo surface of degree 2 and $\s$ is
        the Geiser involution;
\item[6.]   $X$ is a Del Pezzo surface of degree 1 and $\s$ is
        the Bertini involution;
\end{enumerate}
or $n=3$ and $(X, \s)$ is one of the following:
\begin{enumerate}
\item[A1.]   $X$ is a Del Pezzo cubic surface defined by an equation
        of the form $x^3 = F(y,z,w)$ in $\P^3$, and $\s$ is 
        the restriction of the automorphism of $\P^3$ given by
        $(x,y,z,w) \to (\l x,y,z,w)$, where $\l \ne 1$ is a 3rd-root of unity;
\item[A2.]   $X$ is a Del Pezzo sextic surface defined by 
        an equation of the form 
        $z^3 = F(x,y,w)$ in the weighted projective space
        $\P(1,1,2,3)$ with coordinates $(x,y,z,w)$, and $\s$ is 
        the restriction of the automorphism of $\P(1,1,2,3)$ given by
        $(x,y,z,w) \to (x,y,\l z,w)$, where $\l \ne 1$ is a 3rd-root of unity;
\end{enumerate}
or $n=5$ and $(X,\s)$ is one of the following:
\begin{enumerate}
\item[A3.]   $X$ is a Del Pezzo sextic surface defined by 
        an equation of the form 
        $xy^5 = F(x,z,w)$ in the weighted projective space
        $\P(1,1,2,3)$ with coordinates $(x,y,z,w)$, and $\s$ is 
        the restriction of the automorphism of $\P(1,1,2,3)$ given by
        $(x,y,z,w) \to (x,\l y,z,w)$, where $\l \ne 1$ is a 5th-root of unity;
\item[A4.]   $X$ is the Del Pezzo surface $\Bl_{\S}\P^2$, where
        $\S$ is the set of four points in general position, and
        $\s$ is the lift over $X$ of the birational transformation
        of $\P^2$ given, for suitable coordinates of $\P^2$, 
	by $(x,y,z) \to (x(z-y),z(x-y),xz)$.
\end{enumerate}
Moreover, a smooth sextic surface $X$ in $\P(1,1,2,3)$ admits both 
automorphisms $\s_2,\s_3$ such that, for $i=2,3$, the pair 
$(X,\s_i)$ is as in case Ai if and 
only if, in suitable coordinates $(x,y,z,w)$, $X$ is defined by
$x^6 + xy^5 + z^3 + w^2 = 0$.
\end{thmA}

\begin{nt}\label{nt-X_0}
We will denote by $X_0$ the sextic surface in $\P(1,1,2,3)$ defined by the 
equation $x^6 + xy^5 + z^3 + w^2 = 0$.
\end{nt}

\begin{rmk}
Even if we just assumed that the canonical class is not nef, all
surfaces appearing in the classification have Kodaira dimension $- \infty$.
This is expected, since $\Bir(Y) = \Aut(Y)$
for any smooth surface $Y$ with nef canonical class.
\end{rmk}

\begin{thmB}
Let $(X,\s)$ be as in one of cases A1--A4 of Theorem~A.
Denote by $f : X \to X/\s$ the quotient map.
\begin{enumerate}

\item[B1.]      If $(X,\s)$ is as in case A1, then
$X$ is a ``special'' Del Pezzo surface of degree 3, and $X/\s \cong \P^2$.
Moreover, $f$ is defined by the 2-dimensional
linear subsystem of $|-K_X|$ spanned by the orbits of $(-1)$-curves of $X$ 
(see Definition~\ref{orbits-of-$(-1)$-curves}), and is 
totally ramified over a smooth plane cubic curve. 

\item[B2.]      If $(X,\s)$ is as in case A2, then 
$X$ is a ``special'' Del Pezzo surface of degree 1, and
$X/\s \cong \ov \F_3$, the cone in $\P^4$ over a 
rational twisted cubic. Moreover, $f$ is defined by the 
linear subsystem of $|-3K_X|$ spanned by
$3C_x, 2C_x + C_y, C_x + 2C_y, 3C_y, C_w$ (see Notation~\ref{nt-d=1}),
and is totally ramified over the vertex of the cone and the three-canonical
model of a smooth curve of genus 2. 

\item[B3.]      If $(X,\s)$ is as in case A3, then
$X$ is a ``special'' Del Pezzo surface of degree 1, and 
$X \cong X_0$ ($X_0$ is defined in
Notation~\ref{nt-X_0} above) if and only if
$j(C) = 0$ for some (equivalently, for every) smooth $C \in |-K_X|$.
In all cases, $X/\s$ is isomorphic to the sextic hypersurface
of equation $xu = F(x,z,w)$ in the weighted projective space
$\P(1,2,3,5)$ with coordinates $(x,z,w,u)$;
$X/\s$ can be realized by contracting the curve $G' \cup S_0$
of $Z_{22}$ if $X \cong X_0$, and of $Z_{211}$ otherwise
(see Notation~\ref{nt-elliptic-fibration}).
Moreover, $f$ is defined by the 
linear subsystem of $|-3K_X|$ spanned by 
$5C_x, 5C_y, 3C_x + C_z, 2C_x + C_w, C_z + C_w$
(see Notation~\ref{nt-d=1}), and is totally 
ramified over the singular point of $X/\s$
and a smooth elliptic curve.

\item[B4.]      If $(X,\s)$ is as in case A4, then
$X$ is the Del Pezzo surface of degree 5 , 
and $X/\s$ is the image of $Z_{5511}$ under the
contraction of $G_1' \cup G_2' \cup S_0$
(see Notation~\ref{nt-elliptic-fibration}).
Moreover, $f$ is totally ramified over the two singular points of $X/\s$.
\end{enumerate}
\end{thmB}

\begin{rmk}
The ``specialty'' mentioned in cases B1--B3 is characterized
by the constraints given to the equation defining $X$ (see A1--A3). 
It is known that, apart of the Bertini involution, there
are no other automorphisms on general Del Pezzo surfaces 
of degree 1 or 3 (see~\cite{Koitabashi:88}).
The information on the linear systems defining $f$, given for 
cases~B1--B3, will be used to describe the 
birational transformations they induce on $\P^2$. This is not 
needed for case~B4 (one can show that, in this case,
$f$ is defined a linear subsystem of $|-3K_X|$).
\end{rmk}

\subsection{Analogous constructions in higher dimensions.}

Del Pezzo manifolds of dimension $N \geq 3$ are
classified by Fujita~\cite{Fujita:90}.
Del Pezzo manifolds of degree 3 admit projective embeddings as
cubic hypersurfaces, and we find in those which are defined by an equation 
of the form $x_0^3 = F(x_1,\dots,x_{N+1})$ the analogues of case~A1.

If $X$ is a Del Pezzo manifold of degree 2, then
$X$ is a quartic hypersurface in the weighted projective space
$\P(1^{N+1},2)$, and the linear projection
$\P(1^{N+1},2) \rat \P(1^{N+1})$ induces a double covering of $X$
over $\P^N$. This construction generalizes the Geiser involution.

Let now $X$ be a Del Pezzo manifold of degree 1. $X$ is a 
sextic hypersurface in $\P(1^N,2,3)$. The linear projection
$\P(1^N,2,3) \rat \P(1^N,2)$ induces a double covering of $X$
over the cone over Veronese variety $v_2(\P^N)$. 
This is the higher dimensional analogue of the Bertini involution.
If we additionally assume that $X$ is defined, in suitable weighted
coordinates $(x_1,\dots,x_N,z,w)$, by an equation of the 
form $z^3 = F(x_1,\dots,x_N,w)$, then the linear projection
$\P(1^N,2,3) \rat \P(1^N,3)$ induces a triple cyclic covering of $X$
over the cone over $v_3(\P^N)$. This extends case~A2.
Similarly, we can assume that $X$ is defined 
by an equation of the form $x_2 x_1^5 = F(x_2,\dots,x_N,z,w)$. 
Consider the action of order 5 defined on the ring
$\C[x_1,\dots,x_N,z,w]$ by sending $x_1 \to \l x_1$, where
$\l$ is a 5th-root of unity. Then the equation of $X$ is 
invariant with respect to this action, and the
inclusion of the invariant subring of
$\C[x_1,\dots,x_N,z,w]/(x_2 x_1^5 - F)$
determines a degree 5 cyclic covering of $X$ over the sextic
hypersurface of equation $x_2 u = F(x_2,\dots,x_N,z,w)$
in the weighted projective space $\P(1^{N-1},2,3,5)$ of
coordinates $(x_2,\dots,x_N,z,w,u)$. This extends case~A3.

The remaining case to generalize is A4. 
The only Del Pezzo manifolds of degree 5 are the linear sections of 
the Grassmannian variety $\Gr(2,5)$ parametrizing lines in $\P^4$.
The automorphism $\s$ of $X$ defined as in~A4
extends to the higher dimensional Del Pezzo manifolds in the following way.
Let $\s_1$ be a linear automorphism of order 5 of $\C^5$ admitting 
distinct eigenvalues $\l_j = e^{j2 \p i / 5}$ ($j=0,\dots,4$).
Let $e_0,\dots,e_4 \in \C^5$ be the corresponding eigenvectors, and
consider the basis $\{ e_i \wedge e_j \,|\, 0 \leq i < j \leq 4 \}$
for $\C^5 \wedge \C^5$, hence the corresponding Pl\"{u}cker coordinates
$p_{ij}$ of the projective space $\P^9$ in which $\Gr(2,5)$ is embedded.
The automorphism $\s_1$ induces on $\P^9$ the automorphism $\s_2$,
which sends $p_{ij} \to \l_i \l_j p_{ij}$, and the latter
restricts to an automorphism $\s_3$ of $\Gr(2,5)$. 
Consider the five dimensional linear subspace $P \subset \P^9$ given by
$p_{01} - p_{24} = p_{02} - p_{34} = p_{03} - p_{12} = p_{04} - p_{13} = 0$.
Then $X = P \cap \Gr(2,5)$ is a smooth Del Pezzo surface of degree 5, 
$\s_3$ restricts to an automorphism $\s$ of $X$, and
$(X,\s)$ is as in case~A4. In a similar fashion, $\s_3$
restricts to automorphisms on the invariant Del Pezzo manifolds
of intermediate dimensions.

\subsection{Birational transformations of prime order of $\P^2$}

We first recall the definition of three
celebrated birational involutions 
(Examples~\ref{deJonquieres}--\ref{Bertini}) and describe
the constructions of four other birational transformations
of $\P^2$ (Examples~E1--E4). We would like
to point out that, in these examples, not obvious facts will be claimed.
Justifications of such claims are implicitly contained in the proof
of Theorem~E below.

\begin{eg}\label{deJonquieres}
Let $C$ be a curve of degree $d \geq 3$ with an ordinary multiple point $q$
of multiplicity $d-2$. 
The {\it de Jonqui\`eres involution of degree $d$} maps
a general point $p \in \P^2$ to its harmonic conjugate on the line 
$L$ spanned by $p$ and $q$ with respect to the two residual points
of intersection $q',q''$ of $L$ with $C$.
\end{eg}

\begin{eg}\label{Geiser}
Let $\S \subset \P^2$ be a set of 7 points in general position.
The {\it Geiser involution} maps a general point $p \in \P^2$ to the 
ninth base point of the pencil of cubic 
$|\O_{\P^2}(3) \otimes I_{\S} \otimes I_p|$.
\end{eg}

\begin{eg}\label{Bertini}
Let $\S \subset \P^2$ be a set of 8 points in general position.
The {\it Bertini involution} maps a general point $p \in \P^2$ to the 
additional base point of the net of sextics
$|\O_{\P^2}(6) \otimes I_{\S}^2 \otimes I_p|$.
\end{eg}

\begin{egE1}\label{eg-d=3}
Let $\S \subset \P^2$ be a set of 6 points $p_{\a}$ 
such that $X = \Bl_{\S}\P^2$
is as in~A1. Let $\g_1,\dots,\g_6$ 
be the 6 conics passing through all but one point of $\S$,
and $\LL$ denote the set of 15 lines passing through two of the six
points of $\S$. Then $\LL$ splits in the union of subsets
$\LL = \LL' \cup \LL''_1 \cup \LL''_2 \cup \LL''_3$, where 
$\LL' = \{ L_1,\dots, L_6 \}$ and 
$\LL''_{\b} = \{ L_{\b,1},L_{\b,2},L_{\b,3} \}$
($\b = 1,2,3$), such that $L_{\a}$ is tangent to $\g_{\a}$ at $p_{\a}$
(for $\a = 1,\dots,6$) and the three lines $L_{\b,1},L_{\b,2}$ and 
$L_{\b,3}$ meet in one point (for $\b = 1,2,3$).
The set of cubics of the form 
$D_{\a} = \g_{\a} + L_{\a}$ and
$D_{\b} = L_{\b,1} + L_{\b,2} + L_{\b,3}$ spans
a net $W \subset |\O_{\P^2}(3) \otimes I_{\S}|$. 
Imposing any extra general base point to $W$ gives
two additional base points to the system, 
and permutations of such three points define two elements
of order 3 of $\Bir(\P^2)$. 
\end{egE1}

\begin{egE2}\label{eg-d=1-p=3}
Let $\S \subset \P^2$ be a set of 8 points such that $X = \Bl_{\S}\P^2$
is as in~A2. Let 
$\G_w \in |\O_{\P^2}(9) \otimes I_{\S}^3|$ be a curve not contained in
the span of the image of the triple embedding of 
$|\O_{\P^2}(3) \otimes I_{\S}|$ in $|\O_{\P^2}(9) \otimes I_{\S}^3|$,
and let $W$ be the linear subsystem of $|\O_{\P^2}(9) \otimes I_{\S}^3|$ 
spanned by the image of $|\O_{\P^2}(3) \otimes I_{\S}|$ and $\G_w$.
For a suitable choice of $\G_w$, any extra
general base point imposed to $W$ carries with it 
two additional base points, 
and permutations of such three points define two elements
of order 3 of $\Bir(\P^2)$. 
\end{egE2}

\begin{egE3}\label{eg-d=1-p=5}
Let $\S \subset \P^2$ be a set of 8 points such that $X = \Bl_{\S}\P^2$
is as in~A3. There are curves 
$\G_x,\G_y,\G_z,\G_w$ in $|\O_{\P^2}(3a) \otimes I_{\S}^a|$ 
(with $a=1,1,2,3$, respectively) such that, if $W$ is the linear subsystem of
$|\O_{\P^2}(15) \otimes I_{\S}^5|$ spanned by 
$5\G_x, 5\G_y, 3\G_x + \G_z, 2\G_x + \G_w, \G_z + \G_w$,
then for any extra general base point imposed to $W$ the base locus
incorporates four additional base points, and permutations of 
such five points define four elements of order 5 of $\Bir(\P^2)$.
\end{egE3}

\begin{egE4}\label{eg-d=5}
In coordinates $(x,y,z)$ of $\P^2$, let
$\t : (x,y,z) \to (x(z-y),z(x-y),xz)$. Then
$\t$ is an element of order 5 of $\Bir(\P^2)$.
\end{egE4}

\begin{thmE}\label{classification}
Examples~E1--E4 above do define
birational transformations of $\P^2$.
Any element of prime order of $\Bir(\P^2)$ 
is conjugate to one and only one of the birational
transformations described in Examples~\ref{deJonquieres}--E4,
or to an element of $\Aut(\P^2)$. 
Moreover, the transformations defined in these examples, with
the possible exception of Example~E4, are not 
conjugate to elements of $\Aut(\P^2)$.
\end{thmE}

\begin{rmk}
It would be interesting to determine whether Example~E4
is conjugate to an automorphism of $\P^2$.
\end{rmk}

The following theorem gives a
description of the moduli spaces of conjugacy classes
of prime order cyclic subgroups of $\Bir(\P^2)$.
Let $\t \in \Bir(\P^2)$ be an element and $G_{\t} \subset \Bir(\P^2)$ be the 
cyclic subgroup generated by $\t$. 
As in~\cite{BB:99}, associated to any such $\t$
we consider the {\it normalized
fixed curve} $\NFC(\t)$. This is defined as the isomorphism
class of the union of the irrational components of the 
normalization of the curve fixed by $\t$.
Since this is a birational invariant and is the same for
every generator of the group $G_{\t}$, we can define the correspondence 
$$
\NFC : [G_{\t}] \to \NFC(\t)
$$
that associates to any conjugacy class $[G_{\t}]$ of a cyclic subgroup 
of $\Bir(\P^2)$ the normalized fixed curve $\NFC(\t)$ of any generator
$\t$ of any representative $G_{\t}$ of the class.

\begin{thmF}\label{thmF}
The map $\NFC$ naturally establishes one--to--one correspondences between:
\begin{enumerate}
\item[1.]
conjugacy classes $[G_{\t}]$, where $\t$ is a
de Jonqui\`eres involutions of $\P^2$ of degree
$d \geq 3$, and isomorphism classes of (hyper-)elliptic
curves of genus $d-2 \geq 1$;
\item[2.]
conjugacy classes $[G_{\t}]$, where $\t$ is a
Geiser involutions of $\P^2$, and isomorphism classes 
of non-hyperelliptic curves of genus 3;
\item[3.]
conjugacy classes $[G_{\t}]$, where $\t$ is a
Bertini involutions of $\P^2$, and isomorphism classes of 
non-hyperelliptic curves of genus 4 whose canonical
model lies on a singular quadric;
\item[F1.]   
conjugacy classes $[G_{\t}]$, where $\t$ is
as in Example~E1, and isomorphism classes of elliptic curves;
\item[F2.]
conjugacy classes $[G_{\t}]$, where $\t$ is
as in Example~E2, and isomorphism classes of smooth curves of genus 2;
\item[F3.]   conjugacy classes $[G_{\t}]$, where $\t$ is
as in Example~E3, and isomorphism classes of elliptic curves.
\end{enumerate}
\end{thmF}

\begin{rmk}
The four birational transforms $\t,\t^2,\t^3,\t^4$, where $\t$ is
as in Example~E4, form a single conjugacy class 
(see Remark~\ref{conj-d=5} below). As pointed out to us by A.~Beauville, 
linear automorphisms of order $n$, for any given $n < \infty$, form a single 
conjugacy class.
\end{rmk}

Proofs of Theorems E and F are contained in the last section.

\section{Notation and conventions}

We work over an algebraically closed field $k$ of characteristic zero.
We use standard notation accordingly to~\cite{Hartshorne:77},~\cite{BPV:84}
and~\cite{KM:98}.
We will adopt analogous notation as in~\cite{MP:86} 
to denote certain elliptic fibrations over $\P^1$. In particular, 

\begin{nt}\label{nt-elliptic-fibration}
$Z_{22}$ is the (unique) elliptic fibration having one section $S_0$, 
one singular fiber of type $II$, and one of type $II^*$.
$Z_{211}$ is the (unique) elliptic fibration having 
one section $S_0$, two singular fibers of type $I_1$,
and one of type $II^*$.
$Z_{5511}$ is the (unique) elliptic fibration having
five (disjoint) sections $S_0,\dots,S_4$, 
two singular fibers of type $I_5$, and two of type $I_1$.
$Z_{5511}$ is obtained by resolving the base locus of the plane
cubic pencil of equation $y(x-y)(x-z) + \l xz(y-z)=0$.
We fix the following special notation:
if $F_0$ is a fiber of type $II^*$, then we write $F_0 = G \cup G'$
(set-theoretically), where $G$ is the irreducible component 
occurring with multiplicity 5, and $G'$ is the residual curve;
if $F_i$ is a fiber of type $I_5$ of $Z_{5511}$, 
then we denote with $G_i'$ the union of the four components
of $F_i$ which are disjoint from $S_0$ ($i=1,2$).
\end{nt}

\begin{rmk}
To keep this paper more self-contained, 
we opted not to rely on the classification
of extremal rational elliptic fibrations,
which is given in~\cite{MP:86}. We will only use the fact
that the elliptic fibrations
$Z_{22}$, $Z_{211}$ and $Z_{5511}$ are characterized by
the data given in Notation~\ref{nt-elliptic-fibration}
\end{rmk}

\section{On the action of finite groups}

\subsection{Cyclic coverings}

Let $f : X \to Y$ be a cyclic cover,
with $X$ a smooth projective variety of dimension $N$.
Let $R \subset X$ be the set of points fixed by the Galois
action, and set $B = f(R)$. Let
$R = \sum R_k$ and $B = \sum B_k$ be the decompositions
of these cycles with respect to the dimension of their components,
the bottom index standing for the dimension. Actions of finite groups
on smooth varieties can be locally linearized, up to passing to 
completion (see~\cite{BPV:84}, page~85).
This fact implies that $R$ is a smooth cycle and 
$\Sing Y \subseteq (B_{N-2} \cup \cdots \cup B_0)$. 
In particular, $Y$ is smooth if $N=1$, and $Y$ has only isolated singularities
(contained in $B_0$) if $N=2$.

\subsection{Automorphisms on a Del Pezzo surface of degree 1}

If $X$ is a smooth Del Pezzo surface of degree 1, 
its anticanonical ring gives an embedding of $X$
into a weighted projective space $\P = \P(1,1,2,3)$. 
We identify $X$ with its image into this space.
We can write $\P = \Proj S$, where
$S = \C[x,y,z,w]$ is the ring graded by the conditions
$\deg x = \deg y =1$, $\deg z = 2$ and $\deg w =3$. 

Let $\s$ be an automorphism of $X$ of finite order $n$.
Note that $\s$ lifts to 
an automorphism of the sheaf of differentials $\Omega_X$ of $X$.
We deduce that $\s$ is extends to an
automorphism of $\P$. In fact, we can find an automorphism
$\s'$ of $S$, of order $n$, that induces such automorphism of $\P$.
One should observe that the choice of $\s'$ is not unique.

\begin{nt}\label{nt-d=1}
We can assume that the coordinates $(x,y,z,w)$, chosen for $\P$, are 
equivariant, that is, that they define invariant
divisors on $\P$. (Existence of such coordinates
follows from the fact that
linear automorphisms of finite order can be diagonalized.)
We write $\l_x, \l_y, \l_z, \l_w$ for the associated eigenvalues,
so that $\s'$ is given by $(x,y,z,w) \to (\l_x x, \l_y y, \l_z z, \l_w w)$.
The divisors of $\P$ defined by the above coordinates, and 
their restrictions to $X$, will be respectively denoted by
$P_x, P_y, P_z, P_w$ and $C_x, C_y, C_z, C_w$.
Note that the last four divisors
are invariant elements in $|-aK_X|$, for opportune values of $a$.
\end{nt}

To simplify the statement of the following proposition, we rename
$x,y,z,w$ by $x_0,\dots,x_3$, and put $a_i = \deg x_i$.
We write $\l_i$ for the eigenvalue associated to $x_i$. 
We denote by $P_i$ the divisor on $\P$ defined by $x_i$, and 
set $C_i = P_i \cap X$. 

\begin{prop}\label{eigenvalue=1}
Assume that, for some value of $i$, $C_i$ is a fixed divisor of $X$ 
(that is, every point of $C_i$ is fixed).
Then there is an opportune choice of $\s'$
such that $\l_j = 1$ for all $j \ne i$.
\end{prop}

\begin{proof}
Up to renaming the coordinates, we can assume that $i=0$.
We can find a point $p \in C_0$ whose coordinates 
$(0,x_1,x_2,x_3)$ satisfy $x_1 x_2 x_3 \ne 0$.
Since $\s(p) = p$, we have
$\l_1^{1/a_1} = \l_2^{1/a_2} = \l_3^{1/a_3}$.
We can also assume that $a_1 = 1$, and
choose $\s'$ such that $\l_1 = 1$
(this can be done by replacing $\l_i$ with $\l_i\l_1^{-a_i}$).
We deduce then that $\l_2 = \l_3 = 1$, too.
\end{proof}

Clearly, the choice of $\s'$ as in the 
statement of Proposition~\ref{eigenvalue=1} is unique.
We can view this condition as a way of
fixing a distinguished action of $\s$ on $S$.
In fact, if $s \in S$ is the homogeneous element defining $X$,
so that $X = \Proj S/(s)$, then we also get a distinguished
action of $\s$ on $S/(s)$. 
Therefore, when the assumptions of the proposition above
are satisfied, we shall use $\s$ to denote 
these particular automorphisms of $S$ and $S/(s)$. 
Recall that $S/(s)$ is isomorphic to the anticanonical ring of $X$.
It is important to keep in mind that the above ``distinguished'' action 
on $S$ is not necessarily the one induced,
via this isomorphism, by the natural action on the anticanonical ring.

\subsection{Automorphisms of prime order of plane cubic curves}

Let $C$ be an irreducible plane cubic curve 
and $\s$ an automorphism of $C$
of prime order $n \geq 3$. Let $f: C \to C/\s$ be 
the projection on the quotient and $R \subset C$
be the fixed point set.
If $C$ is singular, let $q$ denote the singular point.

\begin{prop}\label{sing-elliptic}
If $C$ is smooth, then either
$R = \emptyset$ and $C/\s$ is an elliptic curve (any value of $n$ may occur),
or $R \ne \emptyset$, $C/\s \cong \P^1$ and $n=3$.
If $C$ is a nodal cubic, then $q$ is the only fixed point
and $C/\s$ is isomorphic to a nodal cubic:
the two tangent directions at $q$ are fixed by $\s$, and the 
corresponding eigenvalues are $\l$ and $\l^{-1}$, where
$\l \ne 1$ is a $n^{th}$-root of unity.
If $C$ is a cuspidal cubic, then there is another
fixed point distinct from $q$ and $C/\s \cong \P^1$.
\end{prop}

\begin{proof}
The first case follows directly by Hurwitz formula. So, assume that
$C$ is singular. Clearly $q$ is a fixed point. 
Take the normalization $\P^1 \to C$.
Then $\s$ lifts to a (non-trivial) automorphism of $\P^1$.
We start considering the case when $q$ is a node of $C$. Let
$q', q'' \in \P^1$ be the two inverse images of $q$. Since $n$ is odd,
$q'$ and $q''$ are necessarily fixed points. Therefore $q$ is the  only fixed
point of $C$ and, locally, the two branches of $C$ passing through $q$
are stabilized by $\s$. By expressing $\s$ in two 
affine charts of $\P^1$, we determine
the two eigenvalues of the action induced on $T_q C$. Then the 
computation of the invariant subring of the local ring 
of $X$ at $q$ shows that $C/\s$ has an ordinary node at $f(q)$.
Now suppose that $q$ is a cuspidal point. By considering the action
on $\P^1$, we see that there is another fixed point beside $q$.
This time the local computation shows that $C/\s$ is smooth at $f(q)$.
\end{proof}

\subsection{Equivariant Mori theory}

Let $X$ be a smooth projective variety
and $G$ be a finite group acting on $X$. Then
$G$ acts on $N^1(X)$ and $N_1(X)$,
and the perfect pairing $(\; \cdot \; ) : N^1(X) \times N_1(X) \to \R$
restricts to a perfect pairing 
$N^1(X)^G \times N_1(X)^G \to \R$.
In particular, $N^1(X)^G$ and $N_1(X)^G$ have the same
dimension, that we shall denote by $\r(X)^G$.

Assume that $\r(X)^G \geq 2$. 
The cone $\CNE(X)^G := \CNE(X) \cap N_1(X)^G$ is called 
{\it $G$-invariant cone of curves} of $X$.
Let $F^G$ be a negative extremal face of $\CNE(X)^G$.
First of all, note that $F^G$ is contained in the boundary of $\CNE(X)$.
Let $F$ be the smallest extremal face of $\CNE(X)$ containing $F^G$.

\begin{prop}\label{equiv-contraction}
$F$ is invariant under the action of $G$, and the extremal contraction
$\cont_F : X \to Y$ of the face $F$ is a $G$-equivariant morphism.
\end{prop}

\begin{proof}
Let $L_F$ be a good-supporting divisor for the face $F$, and
consider the divisor $L  = \sum_{g \in G} gL_F$.
$L$ is still nef, is positive on $\CNE(X) - F$, and vanishes 
along $F^G$. Therefore $L$ is a good-supporting divisor
for some extremal face $F'$ of $\CNE(X)$ with $F^G \subseteq F' \subseteq F$.
By assumption on the dimension of $F$, we conclude that
$F' = F$. Since $L$ is $G$-equivariant, so are $F$ and $\cont_F$.
\end{proof}

\subsection{Resolution of indeterminacy of pairs}

Let $\ff : X' \rat X$ be a birational map between two projective varieties.
If $G$ is a subgroup of $\Bir(X)$ and $G' = \ff^{-1} G \ff$, 
then we say that $(X,G)$ and $(X',G')$ are {\it birationally equivalent pairs}.
If in addition $X'$ is smooth and $G' \subseteq \Aut(X')$, 
then we say that $(X',G')$ is a
{\it resolution of indeterminacy of the pair} $(X,G)$.
We recall the following result (see~\cite{dFE:01} for a proof
and a stronger statement for the case of smooth surfaces).

\begin{thm}\label{resolution}
In the notation above, assume that $G$ is finite. Then
there exists a resolution of indeterminacy $(X', G')$ of the pair $(X,G)$.
\end{thm}

\section{Proofs of Theorems A and B}

Let $X$ be a smooth projective surface whose canonical class
is not nef, and let $\s \in \Aut(X)$ be an element of finite order.

\begin{lem}\label{minimal-pair}
The pair $(X,\s)$ is minimal if and only if for each $(-1)$-curve $E$ of $X$
there exists an integer $k$ such that $E$ and $\s^k E$ intersect properly.
\end{lem}

\begin{proof}
Suppose that there is a $(-1)$-curve $E$ such that, for some $m \geq 1$,
$\s^k E \cap E = \emptyset$ for $1\leq k \leq m-1$ and $\s^m E = E$.
Then the contraction of the disjoint $(-1)$-curves $\s^k E$,
for $k=1,\dots,m-1$, gives an equivariant birational morphism onto
a smooth surface, so $(X,\s)$ is not minimal.
Conversely, let $f : X \to X'$ be a non-trivial birational morphism
of smooth surfaces, $\s' \in \Aut (X')$ be an element 
of finite order, and $\s \in \Aut (X)$ be such that 
$f  \s = \s'  f$. Let $E$ be a $(-1)$-curve contained in
the exceptional locus of $f$ and $m$ the least positive integer such 
that $\s^m E = E$. Note that the curve $C = E + \s E + \cdots + \s^{m-1} E$
is $f$-exceptional, so $C^2 < 0$. 
Suppose that $E$ and $\s^k E$ meet
properly for some positive $k \leq m-1$. 
Then for every $i=0,\dots,m-1$, if 
$j(i)$ denotes the least non-negative integer such that
$i + k \equiv j(i)$ modulo $m$, we have $\s^i E \cdot \s^{j(i)} E \geq 1$.
Note that $j(i) \ne i$. Then we get the contradiction
$$
0 > C^2 = \sum_{i=0}^{m-1} \sum_{j=0}^{m-1}\s^i E \cdot \s^j E \geq 
\sum_{i=0}^{m-1} ((\s^i E)^2 + \s^i E \cdot \s^{j(i)} E) \geq - m + m \geq 0.
$$
\end{proof}

\subsection{The invariant part of the cone of curves}

Assume that $K_X$ is not nef, $\s$ has
prime order $n$, and $(X, \s)$ is a minimal pair.
Recall that $\rho(X)^{\s}$ denotes the rank of $\NS(X)^{\s}$.

\begin{prop}\label{ro=2}           
Suppose that $\r(X)^{\s} \geq 2$.
Then $(X,\s)$ is one of cases~2 and~3 of Theorem~A.
\end{prop}

\begin{proof}
Since $\r(X)^{\s} \geq 2$ and $K_X$ is not nef, we can
find an extremal ray $R$ of $\CNE(X)^{\s}$ contained
in the negative part of $\partial \CNE(X)$. 
Let $F$ be the smallest extremal face of 
$\CNE(X)$ containing $R$. By Proposition~\ref{equiv-contraction},
the contraction of $F$ is $\s$-equivariant. 
We claim that either $F = R$, or $n=2$ and the extremal rays of 
$F$ are generated by $(-1)$-curves $E$ satisfying $E \cdot \s E = 1$.
To see this, suppose that $F \ne R$. Let $[E]$ be a minimal generator of an 
extremal ray of $F$ and consider the curve
$C = E + \s E + \cdots + \s^{n-1}E$.
Note that $[C] \in \partial \CNE(X)$, so $C^2 \leq 0$ 
by~\cite[Lemma~2.5]{Mori:82}. If $n=2$, then 
$$
0 \geq C^2 = E^2 + 2E \cdot \s E + (\s E)^2 = 2(E^2 + E \cdot \s E).
$$
We observe that the possibility $E^2 = E \cdot \s E = 0$ 
can not occur, since $[E]$ and $[\s E]$ 
span different rays. Thus $E^2 = -1$ and $E \cdot \s E = 1$, as claimed.
Consider now the case $n \geq 3$. 
First of all, note that the $\s^i E$ are $(-1)$-curves.
By Lemma~\ref{minimal-pair}, there is a positive integer
$k$ such that $E \. \s^k E \geq 1$. We put $\e = \s^k$.
Note that $\e$ has order $n$, and $(X,\e)$ is a minimal pair.
We deduce that $\e^i E \. \e^j E \ge 0$ for $0 \le i < j \le n-1$,
and in fact $\e^i E \. \e^{i+1} E \ge 1$ for all $i$.
Note also that $C = E + \e E + \dots + \e^{n-1}E$.
Keeping in mind that $n \geq 3$, we get the contradiction
$$
0 \geq C^2 \geq \sum_{i=0}^{n-1}((\e^i E)^2 + 2 \e^i E \cdot \e^{i+1} E) 
\geq -n + 2n > 0.
$$

Now we can conclude the proof of Proposition~\ref{ro=2}. 
First, suppose that $R = F$.
Since $(X,\s)$ is a minimal pair, the $\s$-equivariant contraction 
of $F$ is a $\P^1$-bundle fibration.
If $\s$ acts on each fiber, then there are not
fixed fibers since the set of fixed points of $X$ is smooth and
contains a double section. This gives case~2 of Theorem~A.
Now suppose that $R \subsetneq F$. By the 
above arguments, the $\s$-equivariant contraction 
of $F$ is a conic bundle whose singular fibers are curves
of the form $E \cup \s E$. Let $q$ be the singular point of a singular
fiber. Since the action on $\T_q X$ has
eigenvalues $1$ and $-1$, we can find, locally near
$q$, a fixed section through $q$. This implies that
the automorphism induced on $Y$ is trivial. The same argument used 
for the case of $\s$-invariant $\P^1$-bundles shows that 
$\s$ restricts to an effective action on each fiber.
This gives case 3 of Theorem~A.
\end{proof}

\begin{prop}\label{ro=1,quadric}
Suppose that $X \cong \P^1 \times \P^1$ and $\r(X)^{\sigma} =1$. 
Then $n=2$ and $\s$ swaps the two rulings of $X$.
This is case 4 of Theorem~A.
\end{prop}

\begin{proof}
Since $\r(X)^{\s} = 1$, $n=2$ and $\s$ swaps the two extremal 
rays of $\CNE(X)$.
\end{proof}

\begin{defi}\label{orbits-of-$(-1)$-curves}
Let $X$ be a smooth Del Pezzo surface. 
The set of $(-1)$-curves of $X$ splits in orbits
under the action of $\s$.
We call {\it orbits of $(-1)$-curves} the divisors on $X$ of the form
$D = E + \s E + \cdots + \s^{n-1} E$, with $E$ a $(-1)$-curve of $X$.
\end{defi}

\begin{prop}\label{ro=1}
Suppose that $X \not \cong \P^2, \P^1 \times \P^1$ and $\r(X)^{\s} =1$. 
Then either
$n=2$ and $X$ is a Del Pezzo surface of degree 1 or 2, or
$n=3$ and $X$ is a Del Pezzo surface of degree 1 or 3, or
$n=5$ and $X$ is a Del Pezzo surface of degree 1 or 5.
\end{prop}

\begin{proof}
By pulling back an ample line bundle from $X/\s$, we see that one of
the two generators of $\NS(X)^{\s} \cong \Z$ is ample. 
Since $K_X$ is in $\NS(X)^{\s}$ and is not nef, $-K_X$
is ample, that is, $X$ is a Del Pezzo surface.
Then $X \cong \Bl_{\S}\P^2$ where $\S \subset \P^2$ is a set of $r$ distinct
points in general position, $1 \leq r \leq 8$.
We recall that the degree of $X$ as Del Pezzo surface is $d = K_X^2 = 9 - r$.
Let $D$ be an orbit of $(-1)$-curves of $X$. 
Since $D$ is invariant under the action of $\s$,
$D \in |-a K_X|$ for some positive integer $a$. Then
$n = D \cdot (-K_X) = ad$,
hence $d = 1$ or $n$, since $n$ is prime.
We can conclude by the fact that 
$n$ divides the number of $(-1)$-curves of $X$
(this number is well known as a function of $d$, see
for instance~\cite{Manin:74}). 
\end{proof}

The remaining part of this section is devoted to the analysis
of the six cases presented by Proposition~\ref{ro=1}.

\subsection{The Geiser and Bertini involutions}

Let $\s$ be a biregular involution on a Del Pezzo surface $X$ of
degree 1 or 2, and assume that $\r(X)^{\s} = 1$. Then the action
that $\s$ induces on $\Pic(X)$ is the same 
as the one induced by the Bertini or the Geiser involution,
respectively.
Since their difference induces an automorphism of $\P^2$
which fixes the points of $\S$, they are the same automorphism.

\subsection{Cases A1 and B1}

Let $X$ be a smooth Del Pezzo surface
of degree 3 and $\s$ be an automorphism of $X$ of order 3 such that
$\r(X)^{\s} = 1$. 

\begin{prop}\label{d=3}
Let $X$ and $\s$ as above. Then $(X,\s)$ is as in case~A1, 
$X/\s \cong \P^2$, and the quotient map 
$f : X \to X/\s$ is totally ramified over a smooth
plane cubic. Moreover, $f$ is defined by the linear subsystem of
$|-K_X|$ spanned by the orbits of $(-1)$-curves.
\end{prop}

\begin{proof}
We identify $X$ with its 
anticanonical embedding in $\P^3$. 
Since $\s$ acts on $|-K_X|$, $\s$ is the restriction to $X$
of a linear automorphism of $\P^3$. In particular,
the fixed points set $R$ consists of points, lines
and possibly a smooth plane cubic. Since $\r(X)^{\s} = 1$, 
the eigenvalues of the action of $\s$ on $H^2(X,\Z)$ are 1,
three times $\l$, and three times $\l^2$, with $\l = e^{2\pi i/3}$, 
thus the trace is $-2$. Therefore the trace of $\s$ 
acting on $H^*(X,\Z)$ is $0$. By Lefschetz this is the sum
of the Euler numbers of the fixed components of $\s$. 
We conclude that $R$ is a plane section of $X$, hence
$\s$ fixes a plane in $\P^3$. It follows at once that
$X$ has equation of the form $x^3 = F(y,z,w)$ and $f$
is the projection over the plane $x=0$.
Since the orbits of $(-1)$-curves are mapped to lines of this $\P^2$ 
spanning $|\O_{\P^2}(1)|$, $f$ is defined by the 
claimed linear system.
\end{proof}

\subsection{Cases A2 and B2}

Let $X$ be a Del Pezzo surface of degree 1 and $\s$ be an automorphism
of $X$ of order 3 such that $\r(X)^{\s} = 1$. 
As usual, $f$ will denote the quotient map.
Note that the base point $q$ of $|-K_X|$ is fixed by $\s$.
In particular, $\s$ lifts to an automorphism of $Y = \Bl_qX$
which stabilizes the exceptional divisor $E_q$.
We denote by $g: Y \to Y/\s$ the quotient map.

\begin{lem}\label{trivial}
The action induced by $\s$ on $|-K_X|$ is trivial.
\end{lem}

\begin{proof}
Arguing as in the proof of Proposition~\ref{d=3}, we see that
the trace of $\s$ acting on $H^*(X,\Z)$ is $-1$. If $\s$
does not act trivially on $|-K_X|$ the fixed components 
are points or smooth curves of genus $\le 1$. 
Since these have positive Eulen number, this is impossible. 
\end{proof}

\begin{lem}\label{R}
The ramification locus $R$ of $f$ is the disjoint union of $q$ and $C_z$
(see Notations~\ref{nt-d=1}).
Moreover, $f(q)$ is a singularity of $X/\s$ of type $\frac 1 3 (1,1)$.
\end{lem}

\begin{proof}
It follows from Lemma~\ref{trivial} that $E_q$ is a fixed curve on $Y$. 
Since fixed divisors are smooth, 
we deduce that $q$ is an isolated fixed point of $X$.
$f(q)$ is a singularity 
of type $\frac 1 3 (1,1)$ because $\s$ acts trivially on $\P T_qX$.
Let $R_1$ be the union of the 1-dimensional components of $R$. 
Since $R$ intersects the general $C \in |-K_X|$ in three points, 
one of which is $q$, we have $R_1 \cdot C = 2$. 
Noting that $R_1$ is invariant and $\Pic(X)^{\s} = -K_X \Z$, we 
deduce that $R_1 \in |-2K_X|$. Since $R_1$ intersects properly
every $C \in |-K_X|$, it must be $C_z$. 
By Lefschetz, the Euler number of $R$ is $-1$. Since $C_z$
has Euler number $-2$, we conclude that $q$ is the only 
isolated fixed point.
\end{proof}

\begin{prop}\label{prop-d=1-p=3}
Let $X$ be a smooth Del Pezzo surface of degree 1 and
$\s$ be an element of order 3 in $\Aut(X)$. 
Then $(X,\s)$ is as in case~A2.
\end{prop}

\begin{proof}
We can embed $X$ as a hypersurface of degree 6
in the weighted projective space $\P(1,1,2,3)$. 
We fix coordinates $(x,y,z,w)$ according to Notation~\ref{nt-d=1}. 
Note that $C_z$ is a fixed divisor on $X$. Then
Proposition~\ref{eigenvalue=1} implies that $\s$ 
is the restriction to $X$ of the automorphism of
$\P(1,1,2,3)$ given by $(x,y,z,w) \to (x,y,\l z,w)$,
where $\l \ne 1$ is a $3dr$-root of unity.
In particular, we deduce that
$X$ is defined by an equation of the form $z^3 = F(x,y,w)$.
\end{proof}

We deduce the following

\begin{cor}\label{cor-d=1-p=3}
$f$ is the restriction
of the linear projection of $\P(1,1,2,3)$ from the point $(0,0,1,0)$ 
to $P_z$ (see Notation~\ref{nt-d=1}). In particular, 
$X/\s = P_z \cong \P(1,1,3) \cong \ov \F_3 \subset \P^4$, the 
cone over a rational twisted cubic.
In terms of linear systems, $f$ is defined by the linear subsystem
of $|-3K_X|$ spanned by $3C_x,2C_x + C_y, C_x + 2C_y, 3C_y,C_w$.
\end{cor}

If we embed $X$ in $\P^6$ by $|-3K_X|$, then
the orbits of $(-1)$-curves are hyperplanes sections, the quotient map
$f$ is the restriction to $X$ of a linear projection $\p : \P^6 \rat \P^4$, 
and $X/\s = \ov \F_3 \subset \P^4$.
Let $\f : \F_3  \to \ov \F_3$ be the resolution of the 
singularity $p \in \ov \F_3$, and $E = \f^{-1}(p)$ be the exceptional curve.
Then $Y/\s = \F_3$, $g(E_q) = E$ and $\f g = f \ff$. 
Let $R'$ and $B'$ be the ramification
and branch loci of $g$. Note that they are divisors containing
respectively $E_q$ and $E$, and that $g^* B' = 3 R'$.
Using $2 R' = K_Y - g^* K_{\F_3}$, one sees that 
$B' \sim 2 \f^* \O_{\ov \F_3}(1) + E$. Therefore
$f$ is branched along the vertex $p$ of $\ov \F_3$
and a (smooth) curve $B_1 \in |\O_{\ov \F_3}(2)|$. $B_1$ is a sextic and
$3K_{B_1} \cong \O_{\P^4}(1)|_{B_1}$ by adjunction,
so $B_1$ has genus 2. Riemann-Roch formula implies that the
embedding of $B_1$ in $\P^4$ is given by the complete 
linear system $|K_{B_1}|$. This concludes the proof of~B2.

\subsection{Cases A3 and B3}

Let $X$ be a Del Pezzo surface of degree 1 and $\s$ be an automorphism
of $X$ of order 5.

\begin{lem}\label{effective-d=1-p=5}
$\s$ acts effectively on $|-K_X|$, whose invariant curves
are a rational cuspidal curve $C_x$ and an elliptic curve $C_y$
(see Notations~\ref{nt-d=1}).
The ramification locus of $f : X \to X/\s$ 
is the disjoint union of $C_y$ and the cuspidal point 
$q_1$ of $C_x$. The quotient surface $X/\s$ has
a unique singularity at the point $f(q_1)$.
\end{lem}

\begin{proof}
Since the base point $q$ of $|-K_X|$ is fixed, every invariant
but not fixed $C \in |-K_X|$ must be a cuspidal curve
by Proposition~\ref{sing-elliptic}. 
We deduce that $\s$ acts effectively on $|-K_X|$, so that
$C_x$ and $C_y$ are the only two invariant members of $|-K_X|$. 
By computing the Euler number of $X - (C_x \cup C_y)$, 
we deduce that one of the two curves, 
say $C_y$, is fixed, while the other one, $C_x$, is a rational cuspidal curve.
Note that the cuspidal point $q_1$ of $C_x$ is the only isolated fixed
point of $X$, so $f(q_1)$ is the only singularity of $X/\s$.
\end{proof}

\begin{prop}\label{prop-d=1-p=5}
Let $X$ be a smooth Del Pezzo surface of degree 1 and
$\s$ be an element of order 5 in $\Aut(X)$. 
Then $(X,\s)$ is as in case~A3.
\end{prop}
 
\begin{proof}
The proof is analogous to that of Proposition~\ref{prop-d=1-p=5}. 
This time, in the embedding of $X$ in $\P(1,1,2,3)$,
$\s$ is the restriction of an automorphism acting effectively only
on the $y$ coordinate, and $X$ has equation of the form $xy^5 = F(x,z,w)$.
\end{proof}

\begin{prop}\label{quotient-d=1-p=5}
$X/\s$ can be identified with the sextic hypersurface of 
equation $xu=F(x,z,w)$ in the weighted projective space
$\P(1,2,3,5)$ with coordinates $(x,z,w,u)$.
Let $p_1 = (0,1,0,0) \in \P(1,2,3,5)$. Then $p_1 = f(q_1)$ and
is a singularity of $X/\s$ of type $\frac 15(1,4)$.
In terms of linear systems, $f$ is defined by the 
linear subsystem of $|-5K_X|$ spanned by
$5C_x, 5C_y,3C_x+C_z, 2C_x+C_w, C_z+C_w$.
\end{prop}

\begin{proof}
Since the relation $xy^5-F=0$ involves only the 5th power of $y$,
the ring of invariants $T^{\s}$ of $T = \C[x,y,z,w]/(xy^5-F)$
is generated, over $\C$, by the classes of $x,y^5,z,w$. This yields the
projective description of $X/\s$, since $T^{\s} \cong \C[x,z,w][u]$ 
with $u = F(x,z,w) / x$. Clearly $p_1 = f(q_1)$. The 
equation of $X/\s$ near $p_1$ is approximated by the linear equation 
$x=0$ (of weight 1),
so the nature of the point $p_1 \in X$ is the same as the one
of the point $(0,0,1)$ of $\P(2,3,5)$. This is a singularity of type 
$\frac 15(1,4)$. The last statement of the proposition
follows from the fact that $T_5^{\s}$ 
generates the ring $\oplus_{m \geq 0}T_{5m}^{\s}$.
\end{proof}

We can also follow a more geometric approach to study the quotient
$X \to X/\s$, by understanding explicitly the system $|-K_X|$
in terms of the embedding of $X$ in $\P(1,1,2,3)$.
For $t \in \C$, we put $C_t = C_y - t C_x$.
Note that $C_x$ is given by $F(0,z,w)=0$ in $P_x$, and 
$C_t$ by $F(x,z,w)=t^5 x^6$ 
in $P_t$. Here $P_t \subset \P(1,1,2,3)$ is the subspace defined
by $y = t x$. If $L = P_x \cap P_y$, then 
$C_x \cap L = C_y \cap L = q$. We see that
$q \ne (0,0,1,0),(0,0,0,1)$ 
since $H^0(-2K_X)$ and $H^0(-3K_X)$ are globally generated.
Moreover, we observe that $(0,1,0,0)$ is the cuspidal point $q_1$ of $C_x$.

Note that either $j(C_t) \ne 0$ or $j(C_t) = 0$ for all smooth $C_t$.
Indeed the restrictions of $x,z,w$ to $C_t$ define
divisors on $C_t$ which are linearly equivalent to $q, 2q, 3q$,
respectively. 
Thus, by embedding $C_t$ in $\P^2$ via $|\O_{C_t}(3q)|$,
$C_t$ is defined, in an affine chart, by $F(1,z,w)=t^5$. 
Writing this equation in Weierstrass normal form
$(w')^2 = (z')^3 + A z' + B$, we see that $A$ is independent of $t$. 

Let $\p : \P(1,1,2,3) \rat P_y$ be the linear projection from the 
point $q_1 = (0,1,0,0)$. Its restriction to $X$ contracts $C_x$ 
to $q$ and coincides with $f$ outside $C_x$, since
$\p|_{(X - C_x)}$ is finite of degree 5 and $\s$ acts on the
fibers of $\p$. We resolve the indeterminacy of $\p$ by taking the 
weighted blowup of $\P(1,1,2,3)$ at $q_1 = (0,1,0,0)$, with weights 
$1,2,3$. This restricts to the weighted 
blowup $g : \Bl^w_{q_1}X \to X$ of $X$ at $q_1$ with weights $2,3$.
The latter gives a resolution $f_1 : \Bl^w_{q_1}X \to P_y$ of $\p|_X$. 
Let $E_1$ be the exceptional divisor of $g$ 
and $C_x', C_y'$ be the strict transforms of $C_x, C_y$.
Next we take the blowup $\Bl_q P_y$ of $P_y$ at $q$. Let
$E$ and $L'$ be the corresponding exceptional divisor 
and the strict transform of $L$.
Then $f_1$ factors through a morphism $f_2 : \Bl^w_{q_1}X \to \Bl_q P_y$
and $\Bl_q P_y \to P_y$, and, moreover, $f_2(C'_x) = E$ and $f_2(E_1) = L'$.
Indeed $f_1^{-1}(q) = C_x'$, and this curve intersects $E_1 \cong \P(2,3)$ 
at a point distinct from $(1,0)$ and $(0,1)$.
So $\Bl^w_{q_1}X$ is smooth along this curve, and we can apply
the universal property of the blowup to $f_1$. 
Moreover, $f_2(C_x') = E$, since $C_x'$ is the only
curve contracted by $f_1$, and, by construction,
$f_1(E_1) = P_x \cap P_y = L$, so $f_2(E_1) = L'$.

\begin{prop}
There is a morphism $g': \Bl_q P_y \to X/\s$ such that $f g = g' f_2$.
Moreover, $X/\s = \n(Z)$, where
$Z \cong Z_{22}$ ($\cong Z_{211}$) if $j(C) = 0$ ($\ne 0$, respectively)
for some (equivalently, for every) smooth $C \in |-K_X|$, and 
$\n$ is the contraction of $G' \sqcup S_0$ to $p \sqcup p_1$
(see Notation~\ref{nt-elliptic-fibration}).
\end{prop}

\begin{proof}
Let $h : \~ P_y \to P_y$ be the minimal resolution of the two singularities
of $P_y$. Note that $6L$ is a Cartier divisor, and
$h^*(6L) = 6 \~ L + 3F + 4H + 2H'$.
Here $\~L = h_*^{-1} L$, $F = h^{-1}((0,1,0))$ is a $(-2)$-curve,
and $H \cup  H' = h^{-1}((0,0,1))$ is a chain of two $(-2)$-curves.
We compute $\~ L^2 = -1$ by $L^2 = 1/6$.
Let $\~ q = h^{-1}(q)$. Then $h$ lifts to a morphism
$h' : \Bl_{\~ q} \~ P_y \to \Bl_q P_y$, and
the exceptional divisor $\~ E$ of $\Bl_{\~ q} \~ P_y$ 
is mapped isomorphically to $E$.
The strict transform, over $\Bl_q \~ P_y$, of the curve 
$\~ L \cup F \cup H \cup H'$ 
is a chain of four $(-2)$-curves. Its contraction
$\~ g : \Bl_{\~ q} \~ P_y \to X'$ factors through
$h'$ and a morphism
$g' : \Bl_q P_y \to X'$ which contracts $L'$ (generating
a singularity of type $\frac 1 5 (1,4)$) and is an isomorphism outside $L'$.
This gives $g' f_2 g^{-1} = f$, hence $X' = X/\s$.

Let $D_t = (\p|_X)_*C_t$, for $t \in \C$. 
The minimal smooth resolution $\n : Z \to X/\s$ of the base locus
of $f_*|-K_X|$ (which is supported at $p$) factors as $\n = \~ \n \~ g$.
$Z$ is an elliptic fibration with only one section $S_0$ 
(the exceptional divisor of the last monoidal transformation) and having a 
singular fiber $F_0$ of type $II^*$, so it is either $Z_{22}$
or $Z_{211}$ accordingly to the $j$-invariant of the elements 
in $|-K_X|$. We conclude by observing that $\n$ is the
contraction of $G' \sqcup S_0$, where $G'$ is 
as in Notation~\ref{nt-elliptic-fibration}.
\end{proof}

To conclude the proof of~B3, we observe that 
the branch divisor of $f$ is the proper transform of $C_y \subset P_y$, 
so it is an elliptic curve.

\subsection{Cases A4 and B4}

Let $X$ be the Del Pezzo surface of degree 5 and $\s$ be an automorphism
of $X$ of order 5. Although it is known that $\Aut(X) \cong S_5$,
the symmetric group on five letters, we will explicitly tread
also this case. In a suitable coordinate system,
X is the blowup of $\P^2$ along the set $\S$ of four points $p_i$ 
in general position.
Let $L_{i,j}$ be the line passing through the points $p_i$ and $p_j$,
$L'_{i,j}$ be the proper transform of $L_{i,j}$ over $X$,
$E_i$ be the exceptional divisor over the point $p_i$.
The set $\{L'_{i,j},E_i\}_{i,j}$ is the set of $(-1)$-curves of $X$.
One can check that none of the $(-1)$-curves can be invariant
and the five components of each $\s$-orbit of $(-1)$-curves
are configured into a pentagon. We can assume, without lost of generality,
that the two orbits are
$D_1 = L'_{1,2} + E_1 + L'_{1,4} + L'_{2,3} + E_2$ and 
$D_2 = L'_{3,4} + L'_{1,3} + E_4 + E_3 + L'_{2,4}$.
The five points of intersection of $D_1$ and $D_2$
establish (in a obvious way) a one--to--one
correspondence between the components of $D_1$ and $D_2$.
Fix coordinates $(x,y,z)$ in $\P^2$ such that
$p_1 = (1,0,0)$, $p_2 = (0,1,0)$, $p_3 = (0,0,1)$ and $p_4 = (1,1,1)$, and
let $\t$ be the order 5 Cremona transformation defined by
$$
(x,y,z) \to (xz, x(z - y), z(x - y)).
$$

\begin{prop}\label{d=5}
Let $X$ be the Del Pezzo surface of degree 5 and 
$\s$ be an automorphism of $X$ of order 5. In the notation above, 
after suitably reordering the points $p_i$,
$\s$ is the lift on $X$ of $\t$. This is case~A4.
\end{prop}

\begin{proof}
By means of the correspondence described above,
an order 5 automorphism of $X$ is uniquely determined by the action 
induced on the orbit of one $(-1)$-curve. We deduce that
$\t' = \s^m$ for some $m$. Consider the elements $\e,\f,\ff \in \Aut(\P^2)$
determined by the tree permutations $(p_1p_4p_2p_3)$, $(p_1p_3p_2p_4)$
and $(p_1p_2)(p_3p_4)$ on $\S$, and let
$\e',\f',\ff' \in \Aut(X)$ be their lifts over $X$. Then one can check that 
$\s,\s^2,\s^3,\s^4$ are conjugated one into the other by elements
in $\{ \e',\f',\ff' \}$.
\end{proof}

\begin{rmk}\label{conj-d=5}
$\t,\t^2,\t^3,\t^4$ are conjugated one into the other  
by elements in $\{ \e,\f,\ff \}$.
\end{rmk}

By blowing up $X$ along $D_1 \cap D_2$, we obtain an elliptic 
fibration $Y \cong Z_{5511}$.
The two fibers of type $I_5$ are the strict transforms $D_1'$ and $D_2'$ 
of $D_1$ and $D_2$. Let $F_1$ and $F_2$ denote the other two singular 
fibers. Note that $\s$ acts fiberwise on $Y$.
By Euler number computation and Proposition~\ref{sing-elliptic},
we can see that the elliptic fibration is $\s$-invariant and
the nodes $y_1,y_2$ of the two fibers $F_1, F_2$ 
are the only fixed points of $Y$.
Let $g : Y \to Y/\s$ be the quotient map and 
$\n : Z \to Y/\s$ be the minimal resolution of singularities. 

\begin{prop}\label{prop2-d=5-p=5}
$Y/\s$ is an elliptic fibration over $\P^1$ having exactly 
four singular fibers of type $I_1$ and two singular
points of type $\frac 1 5 (1,4)$, and $Z \cong Z_{5511}$.
\end{prop}

\begin{proof}
The elliptic fibration of $Y$ induces
a fibration on the quotient, and each $D_i'$ is mapped to a 
nodal curve. Proposition~\ref{sing-elliptic}, applied to the 
irreducible fibers of $Y$, yield the first part of the proposition.
By resolving the two singularities of $Y/\s'$, we obtain $Z \cong Z_{5511}$.
\end{proof}

If $G_i = h^{-1}_* g(D'_i)$ and $G_i' = h^{-1}g(y_i) \setminus G_i$,
we deduce the following

\begin{cor}
$X/\s \cong \n(Z_{5511})$, where $\n$ is
the contraction of the cycle $S_0 \sqcup G_1' \sqcup G_2'$.
In particular, $X/\s$ has two singularities of type $\frac 1 5 (1,4)$
and $f: X \to X/\s$ is ramified over these two points.
\end{cor}

\subsection{The last statement of Theorem A}

Let $X \subset \P(1,1,2,3)$ be a smooth sextic surface
and, for $i=2,3$, assume there is a $\s_i \in \Aut(X)$
such that $(X,\s_i)$ is as in case~Ai.
The quotient map $f_i : X \to X/\s_i$ 
is the restriction to $X$ of a linear projection 
$\p_{q_i} : \P(1,1,2,3) \rat P_i$, where
$P_i \subset \P(1,1,2,3)$ is a suitable subspace. We can assume that
$P_i \cap X$ is the ramification divisor of $f_i$. 
Then we can choose coordinates
$(x,y,z,w)$ in $\P(1,1,2,3)$ such that
$q_2 = (0,0,1,0)$, $P_2 = \{ z = 0 \}$, $q_3 = (0,1,0,0)$ and 
$P_3 = \{ y = 0 \}$. 
In this coordinate system, $X$ is defined by an equation of the form
$xy^5 + z^3 + G(x,w) = 0$.
Using that $X$ is smooth, one can check that this equation
is reduced to $xy^5 + z^3 + w^2 + x^6 = 0$ by an opportune 
change of coordinates.
The proofs of Theorems A and B are now complete.

\section{Building back the covering}

The following proposition gives 
the construction of $(X,\s)$ (for cases~B1--B3) 
starting from the quotient $X/\s$.

\begin{propC}
\begin{enumerate}
\item[C1.]      Let $B$ be a smooth cubic of $\P^2$. Then there is a 
        triple cyclic cover $f: X \to \P^2$ branched along $B$, and
        $(X,\s)$ is as in case A1 for any generator $\s$
        of the Galois group of the covering.
\item[C2.]      Let $\ov \F_3 \subset \P^4$ be a cone over a rational 
        twisted cubic, and let $B_1 \subset \ov \F_3$ be a smooth curve of 
        genus 2, cut on the cone by a quadric hypersurface of $\P^4$.
        Then there is a triple cyclic 
        cover $f: X \to \ov \F_3$ branched along $B_1$ and the vertex $p$
        of the cone, and 
        $(X,\s)$ is as in case A2 for any generator $\s$
        of the Galois group of the covering.
\item[C3.]      Let $Z = Z_{22}$ or $Z_{211}$, and $\n : Z \to Y$
        be the contraction of $G' \cup S_0$. 
        Let $B_1$ be a smooth member of $|-K_Y|$.
        Then there is a cyclic cover $f : X \to Y$ of degree 5 branched
        along $B_1$ and the singular point $p_1$ of $Y$, and
        $(X,\s)$ is as in case A3 for any generator $\s$
        of the Galois group of the covering. In particular
        $X \cong X_0$ if and only if $Z = Z_{22}$.
\end{enumerate}
\end{propC}

\begin{proof}
Consider case~C1. The section $s \in H^0(\O_{\P^2}(3))$
vanishing along $B$ determines a triple cyclic cover
$f : X \to \P^2$, where
$X = \Spec (\O_{\P^2} \oplus \O_{\P^2}(-1) \oplus \O_{\P^2}(-2))$.
$X$ is smooth and $f$ is totally ramified over $B$. Since
$-K_X = f^* \O_{\P^2}(1)$, it is ample and 
has self intersection $K_X^2 = 3$ by projection formula.
To conclude, let $\s$ be a generator of the Galois group of $X \to \P^2$.
Then $\r(X)^{\s} = \r(\P^2) = 1$. 

For case~C2, consider the blowup $\f : \F_3 \to \ov \F_3$ of 
$\ov \F_3$ at the vertex $p$.
Let $E$ denote the exceptional divisor. Let $H = \f^* \O_{\ov \F_3}(1)$ and
$\~ B = \f^{-1}(B_1 \cup \{ p \})$. Then
$\~ B \in |2H + E| = |3(H - F)|$.
Set $\LL = \O_{\F_3}(H - F)$. The section 
$s \in H^0(\LL^3)$ defining $B_1$ determines a triple cyclic cover
$\~ f : \~ X \to \F_3$ (totally) ramified over $\~ B$, where
$\~X = \Spec (\O_{\F_3} \oplus \LL^{-1} \oplus \LL^{-2})$.
$\~ X$ is smooth since $\F_3$ and $\~ B$ are smooth.
Set $E' := \~ f^{-1}(E)$. Then $E' = (1/3) \~ f^*(E)$ and
$\~ f_* E' = E$. By projection formula,
$(E')^2 = -1$, so there is a morphism $\ff : \~ X \to X$ 
contracting $E'$ to a smooth point. Then we can find a
morphism $f : X \to \ov \F_3$ such that $f \ff = \f \~f$.
By construction, $f$ is a degree 3 cyclic cover 
branched along $B$. Since $-K_X = (1/3) f^* H_{\ov \F_3}$, $-K_X$ is ample.
Moreover, $H_{\ov \F_3}^2 = 3$ yields $K_X^2 = 1$ by projection formula. 
To conclude, let $\s$ be one of the two generators of 
the Galois group of $X \to \ov \F_3$. Then $\r(X)^{\s} = \r(\ov \F_3) = 1$.

For case~C3, set 
$C = \n (G)$ ($G$ is the irreducible component of the fiber $F_0$
occurring with multiplicity 5, and is the only component not contracted
by $\n$). Then $B_1$ is linearly equivalent to $5C$.
The section $s \in H^0(\O_{Y}(-5C))$
vanishing along $B_1$ determines a degree 5 cyclic covering $f : X \to Y$,
where $X$ is the normalization of
$\Spec (\oplus_{m=0}^4 \O_{Y}(-mC))$.
This covering is \'etale outside $B_1 \cup \{ p_1 \}$ and
totally ramified over $B_1$. A local computation over the point $p_1$
shows that $X$ is smooth and $f$ is totally ramified over $p_1$.
We have $-K_X = f^*C$.
Note that $C^2 = 1/5$, so $C^2$ is ample by Kleiman's criterion.
Hence $-K_X$ is ample as well
and, by projection formula, $K_X^2 = 1$. To conclude, let $\s$ be
any generator of the Galois group of $f$. Then
$\r(X)^{\s} = \r(Y) = 1$.
\end{proof}

\section{Counting the number of automorphisms}

Bertini and Geiser involutions are unique on a 
given Del Pezzo surface of degree 1 or 2.
Here we consider the same question for the  
automorphisms described in cases~A1--A4 of Theorem~A.

\begin{propD}
\begin{enumerate}
\item[D1.]   If $X$ is as in~A1,
        there are exactly eight distinct automorphisms
        as in~A1 if $X$ is the Fermat cubic (i.e., if $X$ is defined by
        $x^3 + y^3 + z^3 + w^3 = 0$), and exactly two otherwise.
\item[D2.]   If $X$ is as in~A2,
        it has exactly two distinct automorphisms as in~A2.
\item[D3.]   If $X$ is as in~A3,
        it has exactly four distinct automorphisms as in~A3.
\item[D4.]   If $X$ is as in~A4,
        it has exactly 24 distinct automorphisms as in~A4.
\end{enumerate}
\end{propD}

\begin{proof}
Case~D1 of the proposition is well known
(see for instance~\cite{Segre:42}, page 129). 
Assume then that $X$ is a smooth sextic surface in $\P(1,1,2,3)$.
D2 follows simply by the fact that there is only one
linear projection from $\P(1,1,2,3)$ onto $\P(1,1,3)$.
Concerning case D3, suppose there are two automorphisms 
$\s_1,\s_2 \in \Aut(X)$ of order 5 corresponding to two distinct linear
projections $\p_{q_i} : \P(1,1,2,3) \rat P_i \cong \P(1,2,3)$.
Then we can fix coordinates $(x,y,z,w)$ in $\P(1,1,2,3)$
such that $q_1 = (1,0,0,0)$, $P_1 = \{ x = 0 \}$, $q_2 = (a,1,0,0)$ and 
$P_2 = \{ bx + cy = 0 \}$
(we assume that $P_i \cap X = R_i$, the ramification divisor
of $f_i$). But one can see, by considering how it
reflects on the equation of $X$, that this situation is impossible. 

Lastly, let $X = \Bl_{\S}\P^2$ be the Del Pezzo surface of degree 5.
Let $\LL$ be the set of lines passing through pairs of points of $\S$.
Any splitting of $\LL$ into two ``triangles'' (i.e. into two
sets of three lines with no common points)
determines a decomposition of the set of $(-1)$-curves of $X$
in two ``pentagons''. To 
each distinct splitting of $\LL$ as above there correspond 
four automorphisms of order 5, and conversely.
We conclude counting six possible distinct ways of splitting $\LL$.
\end{proof}

\section{Proofs of Theorems E and F}

The results proved in the previous sections are finally
applied to prove the classification of birational transformations of prime 
order of $\P^2$ (Theorem~E) and of their
moduli spaces (Theorem~F). Part of the
the arguments used in the proofs are taken from~\cite{BB:99}.

\subsection{Proof of Theorem E}

Let $\t$ be a birational transformation of $\P^2$ of prime order $p$.
By Theorem~\ref{resolution}, there is a resolution $(X,\s)$ 
of the pair $(\P^2,\t)$. We can assume that
$(X,\s)$ is a minimal pair. Note that $K_X$ is not nef,
thus we can apply to $(X,\s)$ the results
stated in Theorem~A of Chapter~5.

In case 2 of Theorem~A, $X$ is isomorphic to an Hirzebruch surface $\F_e$
for some $e \geq 0$. We perform elementary transformations to
reduce $(X,\s)$ to $(\F_1, \s')$, with $\s' \in \Aut(\F_1)$, 
in the following way.
We can find a fixed point of $X$ not contained in the $(-e)$-curve. 
Blowing it up and contracting the proper
transform of the fiber, we obtain $\F_{e-1}$ if $e \geq 2$,
or $\F_1$ if $e=0$. Note that $\s$ induces there an automorphism.
Proceeding in this way, we end up with the desired $(\F_1, \s')$.
Finally, contracting the  $(-1)$-curve of $\F_1$, 
we obtain an automorphism of $\P^2$.

We consider now case 4 of Theorem~A. 
By blowing up a fixed point of $X = \P^1 \times \P^1$
and contracting the proper transforms of the two lines through that point,
we see that $\s$ is birationally equivalent to an automorphism on $\P^2$.

Next, let $(X,\s)$ be as in case 3 of Theorem~A. The fixed divisor
$R \subset X$ is a 
smooth curve passing to the singular point of each reducible fibers 
and intersecting each smooth fiber in two distinct points.
After contracting one component of each reducible fiber of $X$, 
we get on a $\F_e$, where $\s$ induces a birational involution
and $R$ is mapped isomorphically.
Performing elementary transformations centered at general points of 
the image of $R$, 
we eventually get a birational involution $\t'$ on $\F_1$.
After this birational modification, 
$R$ is mapped isomorphically to a curve in $\F_1$.
We will still denote this curve by $R$. Let $E$ and $F$ be
respectively the $(-1)$-curve and a fiber of $\F_1$, and
write $R = 2E + rF$. Adjunction formula yields $r = g+2$, where $g$
is the genus of $R$, thus $E \cdot R = g$.
After further suitable elementary transformations, 
we can lower the multiplicity
of intersection at each point of $R \cap E$ until we get that
$R$ and $E$ meet transversally in $g$ distinct points. At this point we
blow down $E$, obtaining a birational involution of $\P^2$. 
This involution fixed a curve of degree $d = g + 2$ 
with an ordinary multiple point
$q$ of multiplicity $g$ as unique singularity, and lets invariant
the lines through $q$. This is a de Jonqui\`eres involution of degree
$d \geq 2$. To conclude this case, we claim that any
de Jonqui\`eres involution of degree 2 is conjugate to an automorphism.
In analogy with Example~\ref{deJonquieres},
a point $q$ is fixed outside a smooth conic $C$.
Let $T_1$ and $T_2$ be the two lines passing through $p$
and tangent to $C$, $q_i$ be the point of contact of $T_i$ with $C$,
and $L$ be the line spanned by $q_1$ and $q_2$. 
Blowing up $\P^2$ at $q_1$ and $q_2$ and contracting the proper transform
of $L$, we get to $\P^1 \times \P^1$, where the de Jonqui\`eres
involution induces an automorphism $\s'$. We observe that
$\r(\P^1 \times \P^1)^{\s'}=1$, which in particular implies that
$(\P^1 \times \P^1,\s')$ is a minimal pair. By Theorem~A,
$\s'$ is the involution which exchanges the two rulings of $\P^1 \times \P^1$,
so it is birationally equivalent to an automorphism of $\P^2$.

Each case among 5--A4 of Theorem~A is 
clearly birationally equivalent to one of the birational transforms
described in Examples~\ref{Geiser}--E4 (in the same order).
The normalized fixed curve $\NFC(\t)$ 
is given by the isomorphism class of the ramification divisor of 
the cover $X \to X/\s$, and in all cases but A4, 
this is an irrational curve. This shows that
the birational transformations $\t$ described in 
Examples~\ref{deJonquieres}--E3 are not conjugate to elements in $\Aut(\P^2)$. 
Finally, by comparing this invariant together with the order
of the transformation, we conclude that all examples
determine different conjugacy classes.

\subsection{Proof of Theorem F}

Let $(\P^2, \t)$ be one of the pairs listed in 
Examples~\ref{deJonquieres}--E4
and $(X,\s)$ denote its minimal resolution of indeterminacy.
We recall that, excluding the de Jonqui\`eres involutions,
$X$ is a Del Pezzo surface and 
the resolution of $(\P^2,\t)$ is given by 
the blowup of $\P^2$ along $\S$. 

\begin{prop}
The correspondence $\NFC$ is surjective.
\end{prop}

\begin{proof}
Let $C$ be an hyperelliptic curve of genus $g \geq 1$
(with abuse of language, among hyperelliptic curve we
include here also elliptic curves).
Let $g^1_2$ be a pencil of degree 2 on $C$, 
and let $p_1,\dots,p_g \in C$ be $g$ distinct points
such that $p_i + p_j \not \in g^1_2$ for all pairs $i,j$
(for $g=1$ this condition is empty).
The morphism defined by the linear system 
$|g^1_2 + p_1 + \dots + p_g|$ maps $C$ to a plane curve of degree
$g+2$ with an ordinary multiple point $q$ of multiplicity $g$
as unique singularity. This shows that $\NFC$ is 
surjective for the case of de Jonqui\`eres involutions of 
degree $d \geq 3$.

For the remaining cases, since $\NFC(\t)$ is
the isomorphism class of the ramification divisor 
of $f : X \to X/\s$, the claim is clear by Proposition~C except
when $\t$ is as in Example~E2. For this case, 
we have to check that any three-canonical model $C \subset \P^4$ of
a curve of genus 2 lies on a cone over a rational twisted cubic.
Let $W$ be the linear subspace of $|3K_C|$ spanned by the 
image of the triple embedding of $|K_C|$ in $|3K_C|$.
Then the linear projection
$\p : \P H^0(3K_C)^* = \P^4 \rat W^* = \P^3$
maps $C$ two--to--one onto a rational twisted cubic.
The required cone is then obtained by taking the closure of $\p^{-1}\p(C)$.
\end{proof}

\begin{prop}
The correspondence $\NFC$ is injective.
\end{prop}

\begin{proof}
Let $\t$ be a de Jonqui\`eres involution of degree $d \geq 3$ and 
$C \in \P^2$ its fixed curve. Let $\n : \F_1 \to \P^2$ be the blowup
at the singular point $p$ of $C$ and $C_0 = \n^{-1}_* C$.
The fibers of $\F_1$ cut on $C_0$ a $g^1_2$, and 
$C_0 \cap E$ is a set of $g$ distinct points $p_1,\dots,p_g$. Clearly 
$p_i + p_j \not \in g^1_2$ for all pairs $i,j$
and $\n|_{C_0}$ is the morphism defined by 
$|g^1_2 + p_1 + \dots + p_g|$. Then
$\t$ is uniquely determined by $C_0 \cap E$. Let $p_{g+1}$ be
another point of $C$ satisfying $p_i + p_{g+1} \not \in g^1_2$
for all $i=1,\dots,g+1$. After a suitable elementary transformation
we can reduce to the case where $C_0$ is embedded in $\F_1$ in such a way that
$C_0 \cap E = \{ p_2,\dots,p_{g+1} \}$. In this way, by performing
further elementary transformations on $\F_1$, we can produce a birational
link between two any de Jonqui\`eres involutions those
fixed curves have isomorphic normalization $C_0$.

Let now $\t$ be one among Examples~\ref{Geiser}--E4.
In all cases, the isomorphism class of $X$ depends bijectively on
the configuration in $\P^2$ of the points of $\S$ up to linear action
on $\P^2$, and the cover $X \to X/\s$ 
is uniquely determined by its branch divisor up to isomorphism of $X/\s$.
To conclude, it remains to check that the four different
coverings over $\P^2$ of the Fermat cubic of $\P^3$ determine
conjugate cyclic subgroups of $\Bir(\P^2)$.
This is clear, since such coverings are transformed one
into the other by automorphisms of $X$.
\end{proof}

\bibliographystyle{plain}                       
\bibliography{cremona}

\end{document}